\documentclass[10pt,a4paper]{amsart}
\date{\today}

\usepackage{latexsym,amsmath,amsfonts,amscd,amssymb, mathrsfs, slashed, amsthm}
\usepackage{hyperref}
\usepackage[english]{babel}

\theoremstyle{plain}  
\newtheorem{theorem}{Theorem}[section]

\newtheorem*{theorem*}{Theorem}

\theoremstyle{remark}

\newtheorem*{claim*}{Claim}

\numberwithin{equation}{section}

\renewcommand{\geq}{\geqslant}

\newcommand*{\longhookrightarrow}{\ensuremath{\lhook\joinrel\relbar\joinrel\rightarrow}}

\begin{document}
\title{Non-existence of orthogonal complex structures on the round 6-sphere}

\author[A. C. Ferreira]{Ana Cristina Ferreira}
\address{Centro de Matem\'{a}tica \\
Universidade do Minho \\
Campus de Gualtar \\
4710-057 Braga \\
Portugal} \email{anaferreira@math.uminho.pt}

\keywords{6-sphere, complex structure, round metric, twistor space.}
\subjclass[2010]{Primary 53C15; Secondary 53C55, 32L25}

\begin{abstract}
In this short note, we review the well-known result that there is no orthogonal complex structure on $\mathbb{S}^6$ with respect to the round metric. 
\end{abstract}

\maketitle

\section{Introduction}

The volume where this note appears is dedicated to the famous Hopf problem, that is, the question whether there is a complex structure on the 6-sphere. 

Here, we will focus on the round 6-sphere, i.e. $\mathbb{S}^6$ equipped with its notable  metric of constant sectional curvature (unique up to rescaling).  This metric is inherited from the ambient Euclidean metric on $\mathbb{R}^7$ and, under the stereographic projection $\mathbb{S}^6\backslash\{\infty\}\longrightarrow \mathbb{R}^6$, can be written in coordinates $(x_1,x_2,x_3,x_4,x_5,x_6)$ as
$$g = 4 \frac{dx_1^2+dx_2^2+dx_3^2+dx_4^2+dx_5^2+dx_6^2}{(1+x_1^2+x_2^2+x_3^2+x_4^2+x_5^5+x_6^2)^2}$$
(with this scaling factor, its sectional curvature is constant and equal to 1). 

Recall that an almost complex structure $J$ on $\mathbb{S}^6$ is an endomorphism of its tangent bundle $J: T\mathbb{S}^6\longrightarrow T\mathbb{S}^6$ such that its square is minus the identity on each fiber, that is, $J_x^2 = - \mathrm{Id}_x$, for all $x\in \mathbb{S}^6$.  The almost complex structure is said to be orthogonal (with respect to the round metric) if it acts as an isometry of $T\mathbb{S}^6$, more concretely, if $g(JX,JY) = g(X,Y)$, for all vector fields of $T\mathbb{S}^6$. The almost complex structure $J$ is said to be integrable if $\mathbb{S}^6$ becomes a complex manifold such that local charts can be found in which $J$ corresponds to multiplication by $i=\sqrt{-1}$ in $\mathbb{C}^3$. In this case, we say that $J$ is a complex structure and that $(g,J)$ is a Hermitian structure on $\mathbb{S}^6$.

The main objective of this article is to review the proofs, placing them in their historical context, of the following result.

\begin{theorem}\label{theorem}
There is no orthogonal almost complex structure on the round 6-sphere which is integrable. 
\end{theorem}

As far as the author can establish, it was widely believed in the mathematical community that this result had been first proven by Claude LeBrun in 1987, \cite{LeBrun}. The author became aware that this was not the case through a comment of Robert Bryant on {\it Math Overflow}, \cite{MathOverflow}.  The same result had already been established more than thirty years earlier by Andr\'e Blanchard using ideas which anticipated twistors, \cite{Blanchard}.

\bigskip


\section{The article of LeBrun}

The article \cite{LeBrun} consists of a three-page text that presents a concise and ingenious proof of Theorem \ref{theorem}. On the one hand, the proof is elementary in that it does not involve any hard ``machinery'' but, on the other hand, since the arguments are very specific to the case under consideration it cannot be trivially generalized to other situations.  In this section, we will work out this proof backwards and fill in some of the details. 

The proof is set by contradiction. Suppose that there exists an orthogonal complex structure on $\mathbb{S}^6$. The strategy is to exhibit an embedding $\tau: \mathbb{S}^6 \longhookrightarrow M$ where $M$ is a K\"ahler manifold and $\tau$ is a holomorphic map. 

Recall that the De Rham cohomology of spheres is as follows:
$$H_{DR}^k(\mathbb{S}^m) = \left\{\begin{array}{ll}
\mathbb{R}, & \mbox{if}~ k=0,m\\
0, & \mbox{otherwise}.
\end{array} \right.$$ 

By exhibiting the map $\tau$, $\mathbb{S}^6$ would become a complex submanifold of a K\"ahler manifold and, therefore, would also be a K\"ahler manifold. However, for such manifolds the Hermitian form $\Omega(X,Y)=g(JX,Y)$ is closed and its cohomology class $[\Omega]$ is a generator of degree 2. But, since $H_{DR}^2(\mathbb{S}^6)=0$, this cannot happen.

The manifold $M$ considered is $\mathrm{Gr}_3(\mathbb{C}^7)$, the complex Grassmannian of 3-planes in $\mathbb{C}^7$. 
Grassmannians (real or complex) are well-studied examples of manifolds. For $k,n\in\mathbb{N}$, $\mathrm{Gr}_k(\mathbb{C}^n)$ is the set of all complex $k$-dimensional linear subspaces of $\mathbb{C}^n$. For instance, if $k=1$, then $\mathrm{Gr}_1(\mathbb{C}^n) = \mathbb{C}P^{n-1}$. Fixing a Hermitian product on $\mathbb{C}^n$, we can write the Grassmannians as classical homogeneous spaces: 
$$\mathrm{Gr}_k(\mathbb{C}^n) = \frac{U(n)}{U(k)\times U(n-k)}.$$ 
We can readily check the duality $\mathrm{Gr}_k(\mathbb{C}^n) = \mathrm{Gr}_{n-k}(\mathbb{C}^n)$ and that their complex dimension is equal to $k(n-k)$. Complex Grassmannians are K\"ahler manifolds and, as with many properties of such manifolds, the tangent spaces are canonical/tautological objects. More precisely, at each $k$-plane $K$ of $\mathbb{C}^n$ its tangent space is $T_K \mathrm{Gr}_k(\mathbb{C}^n) \simeq \mathrm{Hom}(K, \mathbb{C}^n/K)$, the vector space of linear maps from $K$ to $\mathbb{C}^n/K$.

Let $J$ be any almost complex structure on $\mathbb{S}^6$ and fix $x\in \mathbb{S}^6$. For simplicity of notation, we will write $T$ instead of $T\mathbb{S}^6$ and $T_x$ for the fiber at $x$. Since $J_x: T_x \longrightarrow T_x$ is such that $J_x^2 = -\mathrm{Id}_x$, the eigenvalues of its complexification $J: T\otimes \mathbb{C}\longrightarrow T \otimes\mathbb{C}$ are $+i$ and $-i$. We can define 3-dimensional vector subbundles as follows:
$$\begin{array}{lcl}
T^{1,0} & = & \{v\in T\otimes\mathbb{C}| ~ Jv = iv \},\\
T^{0,1} & = & \{v \in T\otimes\mathbb{C}| ~ Jv =- iv\}.
\end{array}$$ 
Since $T_x \subset \mathbb{R}^7$ then $T_x^{0,1} \subset \mathbb{C}^7$ and we have a tautological map
$$\begin{array}{llcl}
\tau: & \mathbb{S}^6 & \longrightarrow & \mathrm{Gr}_3(\mathbb{C}^7)\\
	& x & \longmapsto & T_x^{0,1}
\end{array}$$
(the reason why $T_x^{0,1}$ is chosen instead of $T_x^{1,0}$ will become clear in the sequel). Let us prove that $\tau$ is an embedding. It suffices to construct a left inverse $\pi: \mathcal{U} \longrightarrow \mathbb{S}^6$ such that $\pi\circ \tau$ is the identity map of $\mathbb{S}^6$ and $\mathcal{U}$ is an open set of $\mathrm{Gr}_3(\mathbb{C}^7)$. For a 3-plane $P$
denote by $\overline{P}$ its complex conjugate. Take $\mathcal{U}$ to be the set of $3$-planes such that $P\cap \overline{P} = \{0\}$, then $\tau(\mathbb{S}^6)\subset \mathcal{U}$ (because $\overline{T_x^{0,1}}= T_x^{1,0}$). If $P\in \mathcal{U}$, then $(P+\overline{P})\cap \mathbb{R}^7$ is a 6-dimensional real vector space, and we can define the map $\tilde{\pi}$ as
$$\begin{array}{llcl}
\tilde{\pi}: & \mathcal{U} & \longrightarrow & \mathbb{R}P^6\\
	& P & \longmapsto & [N_P]
\end{array}$$
where $[N_P]$ is set of the two normal unit vectors of $(P+\overline{P})\cap \mathbb{R}^7$ (notice that $(P+\overline{P})\cap \mathbb{R}^7$ is not canonically oriented) . By recalling that $$T_x \mathbb{S}^6 = \{v\in\mathbb{R}^7: v \perp x \}$$ then $\tilde{\pi} \circ \tau : \mathbb{S}^6 \longrightarrow \mathbb{R}P^6$ lifts to a map $\mathbb{S}^6\longrightarrow\mathbb{S}^6$ which is either $\mathrm{Id}$ or $-\mathrm{Id}$. By taking the first alternative, the claim that $\tau$ is an embedding follows. 

Remark that we did not use the assumption of $J$ being integrable to establish the embedding. We need that hypothesis to derive that $\tau$ is a holomorphic map. 

We will now show that $\overline{\partial}\tau = 0$, that is, the $(0,1)$-part of $d\tau$ vanishes. At a point $x\in \mathbb{S}^6$, we have that  $d\tau$ is given by $$d_x\tau: T_x\otimes \mathbb{C} \longrightarrow T_{\tau(x)}\mathrm{Gr}_3(\mathbb{C}^7) = \mathrm{Hom}(T_x^{0,1}, \mathbb{C}^7/T_x^{0,1}).$$ 
For a vector  $v\in T_x\otimes\mathbb{C}$, $d_x\tau(v)$ can be identified with $D_v \, (\mathrm{mod} ~ T^{0,1})$ where $D$ is the extension to $\mathbb{C}^7$ by complex linearity of the Levi-Civita connection of $\mathbb{R}^7$. To understand why this is so, we need to take a closer look at the isomorphism $T_V \mathrm{Gr}_3(\mathbb{C}^7) \simeq \mathrm{Hom}(V, \mathbb{C}^3/V)$. 
Consider the standard action of $\mathrm{GL}(7,\mathbb{C})$ on $\mathrm{Gr}_3(\mathbb{C}^7)$
$$\begin{array}{ccl}
\mathrm{GL}(7,\mathbb{C}) \times \mathrm{Gr}_3(\mathbb{C}^7)   & \longrightarrow & \mathrm{Gr}_3(\mathbb{C}^7) \\
(A,V) & \longmapsto  & AV
\end{array}.$$
which we differentiate at $(I,V)$ to get a linear map
$$\mathrm{gl}(7,\mathbb{C})\times T_V \mathrm{Gr}_3(\mathbb{C}^7) \longrightarrow T_V \mathrm{Gr}_3(\mathbb{C}^7). $$
Plugging the 0-vector into the second argument we get a surjective  map
$$\Psi: \mathrm{gl}(7,\mathbb{C}) \longrightarrow T_V \mathrm{Gr}_3(\mathbb{C}^7).$$
We can also consider the natural map
$$\varphi: \mathrm{gl}(7,\mathbb{C})\longrightarrow \mathrm{Hom}(V, \mathbb{C}^7/V)$$ in which a linear map is sent to its restriction to $V$ composed with the projection to $\mathbb{C}^7/V$. Note that the isotropy of $V$ of the $\mathrm{GL}(7,\mathbb{C})$-action on $\mathrm{Gr}_3(\mathbb{C}^7)$, the subgroup of automorphisms leaving $V$ invariant, is such that its Lie algebra, the set of endomorphisms of $\mathbb{C}^7$ that send $V$ to itself, is the kernel of $\varphi$. 
Thus, we have a sequence of maps  
$$\mathrm{gl}(7,\mathbb{C})\stackrel{\varphi}{\longrightarrow} \mathrm{Hom}(V,\mathbb{C}^7/V)\stackrel{\psi}{\longrightarrow} T_V\mathrm{Gr}_3(\mathbb{C}^7)$$
such that $\Psi = \psi\circ \varphi$. The map $\psi$ gives an isomorphism  $T_V\mathrm{Gr}_3(\mathbb{C}^7) \simeq \mathrm{Hom}(V,\mathbb{C}^7/V).$ 
Given $A\in \mathrm{gl}(7,\mathbb{C})$ and a curve $B(t)\in \mathrm{GL}(7,\mathbb{C})$ such that $B(0) = I$ and $B'(0)=A$, then
$$\Psi(A) = \left.\frac{d}{dt}\right|_{t=0}(B(t)V)$$
and this is known to be independent of the curve chosen. We could choose, for example, $B(t) = I + tA$ (for small t) and then
$$\Psi(A) = \left.\frac{d}{dt}\right|_{t=0} (\{v+tA(v)\, |\, v\in V \}).$$  
Using $\psi$, we thus identify $\Psi(A)$ with the linear map
$$v \longmapsto B'(0)v\, (\mathrm{mod}\,V), \quad v \in V,$$ in $\mathrm{Hom}(V,\mathbb{C}^7/V)$. 
In other words, what the identification does is to interchange differentiation and the insertion of $V$. But $B'(0)v$ is also 
$\left. \frac{d}{dt}\right|_{t=0}B(t) v$; here $B(t)v$ is a vector field that extends $v$ and at the time $t$ takes values in $B(t)(V)$. 

Now, let $c(t)$ be a curve in $\mathbb{S}^6$ such that $c(0) = x$ and $c'(0) = v \in T_x\otimes \mathbb{C}$. Then $d_x\tau(v) = (\tau\circ c)'(0)$. Take a curve $B(t)$ in $\mathrm{GL}(7,\mathbb{C})$ such that $$B(t)(T^{0,1}_{x}) = T^{0,1}_{c(t)} = (\tau\circ c)(t).$$
 Reasoning as above, 
$$\left.\frac{d}{dt}\right|_{t=0}(\tau\circ c)(t) = \left.\frac{d}{dt}\right|_{t=0} (B(t)(T_x^{0,1}))$$
and, after identification, we obtain the linear map
$$w \longmapsto B'(0)w\, (\mathrm{mod}\, T_x^{0,1}).$$
The curve $B(t)v$ is an extension of $v$ taking values in $T^{0,1}_x$. In flat space, the $t$-derivative coincides with the covariant derivative, thus $d\tau_x$ is the map
$$d\tau_x(v): w \longmapsto D_v(w) \, (\mathrm{mod}\, T_x^{0,1}).$$

Returning to the holomorphicity of $\tau$,  what we need to show, in practical terms, is that $D_V W \in T^{0,1}$ for any two vector fields in $T^{0,1}$. Recall that
$\nabla_X Y = D_X Y + g(X,Y)N$ where $\nabla$ is the Levi-Civita connection of $\mathbb{S}^6$, $N$ is the unit outward-pointing normal vector field of $\mathbb{S}^6$ and $X,Y$ are any two vector fields of $\mathbb{S}^6$. For $V,W \in T^{0,1}$, since $J$ is orthogonal with respect to $g$, then
$$g(V,W) = g(JV, JW) = g(-iV, -iW) = -g(V,W)$$
and therefore, $g(V,W) = 0$. In particular, we have that, for $(0,1)$-fields, $\nabla$ coincides with $D$. 

Since $J$ is integrable, there is a system $\{z^1, z^2, z^3\}$ of local holomorphic coordinates of $\mathbb{S}^6$. Let $X_\alpha$, $\alpha =1,2,3$, be the complex vector fields on $\mathbb{S}^6$ determined by the pairing
$$\langle dz^\alpha, Y \rangle = g(X_\alpha, Y), \quad \forall Y\in T\otimes\mathbb{C}.$$ Thus defined, the set $\{X_1,X_2,X_3\}$ forms a local frame of $T^{0,1}$. Also, for any two such fields, the following equality holds: $$\nabla_{X_\alpha}X_\beta = - \nabla_{X_\beta}X_\alpha.$$  In \cite{LeBrun}, a local computation to show this equation is presented. Here is an index-free exposition. Let $V,W$ be any (complex-valued) vector fields and $f$ a function on $\mathbb{S}^6$. We have that $$(\nabla_V df)(W) = (\nabla_W df)(V)$$ and this is seen as follows. Using the fact that $\nabla$ has no torsion and the definition of the Lie bracket of two vector fields, we get $$(\nabla_V df)(W) = V.(df(W))-df(\nabla_V W) = V.(W.f)-df(\nabla_W V + [V,W]) =$$  $$V.(W.f) - [V,W].f - df(\nabla_W V) = W.(V.f) - df(\nabla_W V) = (\nabla_W df)(V).$$   
Since $\nabla$ is a metric connection and $g(X_\alpha,X_\beta) = 0$, for any vector field $V$
$$g(\nabla_V X_\alpha, X_\beta) + g(X_\alpha, \nabla_V X_\beta) = 0.$$
Also
$$g(\nabla_V X_\alpha, X_\beta) = dz^\beta(\nabla_V X_\alpha) =  - (\nabla_Vdz^\alpha)(X_\beta) = -(\nabla_{X_\beta}dz^\alpha) (V).$$ 
Analogously, $g(\nabla_V X_\beta, X_\alpha) = -(\nabla_{X_\alpha}dz^\beta)(V)$. Given that $V$ is arbitrary, then $\nabla_{X_\alpha} dz^\beta = \nabla_{X_\beta} dz^\alpha$. Thus, since $\nabla$ preserves $g$, $\nabla_{X_\alpha}X_\beta = -\nabla_{X_\beta}X_\alpha$.

Using the equality just established, we have that 
$$[X_\alpha, X_\beta] = \nabla_{X_\alpha}X_\beta - \nabla_{X_\beta}X_\alpha = 2 \nabla_{X_\alpha}X_\beta.$$ 
For $X_\alpha, X_\beta \in T^{0,1}$, it follows from the integrability of $J$, that $[X_\alpha, X_\beta]\in T^{0,1}$ and hence that $\nabla_{X_\alpha} X_\beta \in T^{0,1}$, proving that $\tau$ is holomorphic. This concludes the proof of Theorem \ref{theorem}.

\bigskip

\section{The article of Blanchard}

Blanchard's article appeared in 1953, 34 years before the proof of LeBrun. Remarkably, it contains ideas of twistor theory, well before the use of the word ``twistor'' or the Penrose program, although its roots go back much further to Klein and Lie. For an overview of the history of twistor theory, we recommend Penrose's own account, \cite{Penrose-origins}. 

In the survey article \cite{Salamon}, Simon Salamon presents a proof of Theorem \ref{theorem} in the language of twistors. Although the two proofs are not exactly the same, LeBrun's arguments are specific to the 6-sphere whereas a twistor proof is far more generic, Salamon claims they are formally equivalent. Strikingly, a quick read of Salamon's proof reveals that it is, in fact, the proof of Blanchard written in modern terminology.   As mentioned in the introduction, the geometry community was unaware of Blanchard's article until very recently and, indeed, this article is not cited in Salamon's survey.

Another interesting point about Blanchard's text is that it appeared before the Newlander-Nirenberg theorem, \cite{NewlanderNirenberg}. In fact, the article starts with a manifold $V_{2n}$ equipped with an almost complex structure $J$ and a tensor $T$  given in coordinates by
$$T_{ij}^k = (\partial_i J_j^m - \partial_j J_i^m)J_m^k+J_j^m\partial_m J_i^k - J_i^m\partial_m J_j^k.$$
This is the Nijenhuis tensor $N_J$. Blanchard observes that if $J$ is integrable then $T=0$ and that the vanishing of $T$ guarantees the integrability of $J$ in the real-analytic case. Note that this is a consequence of Frobenius theorem and so it does not require the Newlander-Nirenberg result.  

The proof of Blanchard/Salamon goes as follows. Consider a domain $D$ in $\mathbb{R}^{2n}$. Fix the Euclidean metric and an orientation on $D$. Consider $\mathcal{J}$, the set of all possible almost complex structures on $D$ which are compatible with the metric and the orientation. Viewing an almost complex structure $J$ as a matrix such that $J^2=-\mathrm{Id}$, $\mathrm{SO}(2n)$ acts on $\mathcal{J}$ by conjugation and the stabilizer of a given $J$ 
$$\{h\in \mathrm{SO}(2n):~ hJh^{-1} = J  \}$$
is isomorphic to the unitary group $\mathrm{U}(n)$. Thus $\Gamma_n = \mathrm{SO}(2n)/\mathrm{U}(n)$ parametrizes the set of all compatible almost complex structures on $D$ at a fixed point. $\Gamma_n$ has a complex structure defined by the center of $\mathrm{U}(n)$ which is compatible with the metric induced by the Killing form of $\mathrm{SO}(2n)$. Furthermore, this defines a K\"ahler structure on $\Gamma_n$, \cite{BorelLich}.

Blanchard continues with his Euclidean argument, but we make a digression to explain the twistor construction. Given a $2n$-dimensional oriented Riemannian manifold $M$, let $P \longrightarrow M$ denote the principal $\mathrm{SO}(2n)$-bundle of orthonormal positively-oriented frames of $M$. Using the standard associated bundle construction, we can form the vector bundle
$$\pi: P\times_{\mathrm{SO}(2n)} \Gamma_n \longrightarrow M.$$
The total space $Z^+$ of this new bundle is then the twistor space of $M$. Each fiber $\pi^{-1}(x)$ parametrizes the positively-oriented orthogonal almost complex structures on the vector space $T_x M$ and, by construction, an orthogonal almost complex structure $J$ on an open set $U$ of $M$ determines a local section $s_J: U \longrightarrow \pi^{-1}(U)$ of $Z^+$. 

Two very important claims here are the following: 

\begin{theorem}
$Z^+$ admits an almost complex structure $I$ such that $ds_J$ is complex-linear (i.e. $ds_J\circ J = I \circ ds_J$) if and only if $J$ is integrable.
\end{theorem}
 This result was first advanced for $n=2$ by Penrose in \cite{Penrose} with another proof by Atiyah, Hitchin and Singer in \cite{AtiyahHitchinSinger}, and then by Salamon in \cite{Salamon-LNM} in full generality. Furthermore, we have 

\begin{theorem} 
 $(Z^+, I)$ is a complex manifold if and only if:
 
- $M$ is conformally flat (or, equivalently, the Weyl tensor vanishes) for $n\geq 3$; 

- $M$ is anti-self-dual (or, equivalently, the self-dual-part of the Weyl tensor vanishes) for $n=2$.

\end{theorem} 
This last statement appeared, again for $n=2$, in \cite{Penrose, AtiyahHitchinSinger}, and for $n\geq3$ was established by O'Brian and Rawsnley in \cite{OBrianRawnsley}. Remark that, by combining these last two results, we have, for a conformally flat complex manifold $M$, that the section $s_J: M \longrightarrow Z^+$ is a holomorphic embedding. 

Returning to Blanchard's article, we see that first it is established that the twistor space of the domain $D$ (which is simply $D\times \Gamma_n)$ is a complex manifold and that the map $x \longmapsto (x,J_x)$ is  holomorphic. Then, for a conformally flat space, a cover by open sets and a gluing procedure is used to show that the same principle can be established for such manifolds, constructing a fibered space which is effectively the twistor space. Finally, the assertion that an orthogonal complex structure can be identified with a section of such a bundle and that the base space can be seen as a complex submanifold of the total space is stated. What Blanchard did not prove, and for his purposes there was no need to do so, was the converse in the above two theorems.     

Finally, if we look at $\mathbb{S}^6 = \mathrm{SO}(7)/\mathrm{SO}(6)$ with its standard conformally flat metric (the round metric) we have the fibration
$$\mathrm{SO}(6)/\mathrm{U}(3) \longrightarrow \mathrm{SO}(7)/\mathrm{U}(3) \longrightarrow \mathrm{SO}(7)/\mathrm{SO}(6)$$
so the twistor bundle of $\mathbb{S}^6$ has total space $\mathrm{SO}(7)/\mathrm{U}(3) = \mathrm{SO}(8)/\mathrm{U}(4)=\Gamma_4 $ which again is a K\"ahler manifold. Therefore, we arrive at the same type of contradiction as in LeBrun's proof. 

\bigskip

\section{Concluding remarks}

We can play with the same type of arguments for the other even-dimensional spheres. However, except for $n=1$ and $n=3$, we have no sections of the twistor bundle or, in simpler terms, there are no almost complex structures on $\mathbb{S}^{2n}$, \cite{BorelSerre}. For $n=1$, $\mathbb{S}^2= \mathbb{C}P^1$, and there is no contradiction because $H^2_{DR}(\mathbb{S}^2)=\mathbb{R}$. 

As we can see, the properties of the round metric $g$, especially the fact that it is inherited from the ambient metric of $\mathbb{R}^7$, are heavily used in the arguments above. For a different proof, that uses the curvature properties of $g$, see \cite{SekigawaVanhecke}. Yet, as far as the author could determine, all proofs of Theorem \ref{theorem} get their contradiction from the vanishing of the second de Rham cohomology group. 

Another proof by different methods and its generalization to metrics in a neighborhood of the round metric can be found in  the review by B. Kruglikov in this volume.  

\bigskip

\section*{Acknowledgements} A.C. Ferreira acknowledges the organizers of MAM-1 for the opportunity to give a talk and be a part of the task force that composed this volume. She also wishes to express her many thanks to Oliver Goertsches and Boris Kruglikov   for the clarification of some points and discussion of the topic.   

\bigskip

\bibliographystyle{plain}

\bibliography{mybibfile}

\end{document}